\documentclass[pdflatex,sn-mathphys-num]{sn-jnl}

\usepackage{lineno}
 \usepackage{enumerate}
\usepackage{graphicx}%
\usepackage{multirow}%
\usepackage{amsmath,amssymb,amsfonts}%
\usepackage{amsthm}%
\usepackage{mathrsfs}%
\usepackage[title]{appendix}%
\usepackage{xcolor}%
\usepackage{textcomp}%
\usepackage{manyfoot}%
\usepackage{booktabs}%
\usepackage{algorithm}%
\usepackage{algorithmicx}%
\usepackage{algpseudocode}%
\usepackage{listings}%

\def\f{\longrightarrow}
\def\im{\Longrightarrow}

\def\<{\langle}
\def\>{\rangle}

\def\R{\mathbb{R}}

\def\inte{\textnormal{int}\,}
\def\clo{\textnormal{cl}\,}
\def\epi{\textnormal{epi}\,}
\def\bdry{\textnormal{bdry}\,}
\def\dom{\textnormal{dom}\,}

\def\nvarphi{\bp\nvarphi}

\def\bp{\hspace{-0.08cm}}

\def\far{\textnormal{far}}
\makeatletter
\newcommand*\rel@kern[1]{\kern#1\dimexpr\macc@kerna}
\newcommand*\widebar[1]{%
  \begingroup
  \def\mathaccent##1##2{%
    \rel@kern{0.8}%
    \overline{\rel@kern{-0.8}\macc@nucleus\rel@kern{0.2}}%
    \rel@kern{-0.2}%
  }%
  \macc@depth\@ne
  \let\math@bgroup\@empty \let\math@egroup\macc@set@skewchar
  \mathsurround\z@ \frozen@everymath{\mathgroup\macc@group\relax}%
  \macc@set@skewchar\relax
  \let\mathaccentV\macc@nested@a
  \macc@nested@a\relax111{#1}%
  \endgroup
}
\makeatother

\theoremstyle{thmstyleone}%
\newtheorem{theorem}{Theorem}
\newtheorem{corollary}[theorem]{Corollary}
\newtheorem{lemma}{Lemma}%

\theoremstyle{thmstyletwo}%
\newtheorem{example}{Example}%
\newtheorem{remark}{Remark}%

\theoremstyle{thmstylethree}%

\begin{document}
\title[Spherically Supported Sets with Variable Radii]{Spherically Supported Sets with Variable Radii}


\author*[1]{\fnm{Chadi} \sur{Nour}}\email{cnour@lau.edu.lb}

\affil*[1]{\orgdiv{Department of Computer Science and Mathematics}, \orgname{Lebanese American University}, \orgaddress{\city{Byblos Campus}, \postcode{P.O. Box 36}, \state{Byblos}, \country{Lebanon}}}


\abstract{We study closed sets for which, at every boundary point, there
exists a containing supporting ball whose prescribed radius may vary
from point to point. This extends the notion of spherical support from
constant to variable radii. The variable-radius setting presents a new
difficulty: unlike the constant-radius case, the intersection of the
supporting balls need not be strongly convex with respect to the
prescribed radii. We overcome this difficulty by associating with every
such set a canonical set and an induced radius function determined by
the active supporting balls. For continuous prescribed radii, our main result establishes a geometric
dichotomy: the original set either coincides with the canonical set and
is strongly convex with respect to the prescribed radii, or it is
contained in the boundary of the canonical set, which is strongly
convex with respect to the induced radii. An example shows that the
induced radii cannot, in general, be replaced by the prescribed ones in
the second alternative. We also recover the constant-radius theory and characterize
variable-radius strong convexity in terms of the spherical support
of the boundary.}




\maketitle

\section{Introduction}\label{sec1}

\subsection{Geometric motivation and overview}

Supporting hyperplanes are a fundamental tool for describing the
geometry of convex sets. Replacing hyperplanes by spheres leads to
stronger geometric requirements: a supporting ball must contain the
entire set and touch it at the boundary point under consideration.
The radius of the ball measures the strength of this requirement,
with smaller radii imposing a tighter geometric constraint. Such
descriptions are central to the study of strong convexity, a notion
with connections to optimization, variational analysis, and
approximation theory; see \cite{Azagra,frankowska,Nacry2025,weber}.

Two different requirements arise naturally when the supporting radius
is fixed. One may ask for a containing supporting ball of that radius
at every boundary point, or require such a ball in every outward normal
direction at each boundary point. The first requirement defines spherical support, whereas the second
describes strong convexity. Although the second is more demanding,
the two are closely related. For a fixed supporting radius, spherical
support implies strong convexity whenever the set has nonempty
interior. More generally, a spherically supported set either is
strongly convex or lies in the boundary of a strongly convex set;
see \cite{JCA2026,Pos2026}.

Allowing the supporting radius to vary along the boundary provides
a more flexible description of the geometry. Variable-radius strong
convexity was studied in \cite{Nacry2025}. The question addressed
here is whether the geometric relationship between spherical support
and strong convexity persists when the radius is prescribed separately
at each boundary point. A difficulty arises because intersecting
supporting balls with different radii need not preserve strong
convexity with respect to the prescribed radii.

We answer this question by constructing a canonical containing set
as the intersection of the prescribed supporting balls. When the prescribed radius varies continuously,
the original set either coincides with this set and is strongly convex
with respect to the prescribed radii, or lies in its boundary. In the
latter case, the canonical set is strongly convex with respect to
induced radii, whose necessity is illustrated by an example. We also
characterize variable-radius strong convexity through spherical
support of the boundary, extending \cite{Pos2026}. The next subsection provides the precise definitions, the relevant
background, and a detailed presentation of our contributions.


\subsection{Mathematical background and contributions}

A nonempty closed set $S\subset\mathbb R^n$ is said to be
\emph{$r$-strongly convex}, for some $r>0$, if it can be represented as an
intersection of closed balls of radius $r$. Equivalently, at every
boundary point $s$ of $S$ and for \emph{every} unit proximal normal
direction $\zeta$ to $S$ at $s$, one has $S\subset\overline B(s-r\zeta;r)$, where $\overline B(x;\rho)$ denotes the closed ball of radius $\rho$
centered at $x$. This notion admits several other equivalent
geometric and analytic characterizations; see
\cite{Azagra,frankowska,weber} and the references therein. A detailed
historical account can be found in the introduction of
\cite{Nacry2025}.

The notion of an \emph{$r$-spherically supported} set was introduced
in \cite{JCA2026}. A nonempty closed set $S\subset\mathbb R^n$ is
$r$-spherically supported if, at every boundary point $s$, there
\emph{exists} a unit proximal normal direction $\zeta$ to $S$ at $s$
such that $S\subset\overline B(s-r\zeta;r)$. Thus, strong convexity requires this inclusion in every unit proximal
normal direction, whereas spherical support requires it in at least
one such direction at each boundary point. Every $r$-strongly convex
set is therefore $r$-spherically supported. The converse fails in
general, as demonstrated in \cite{JCA2026} by taking $S$ to be the
unit circle in the plane. Nevertheless, Nour and Takche proved in
\cite[Theorem~1.1]{JCA2026} that a nonempty closed set $S$ not reduced
to a singleton is $r$-strongly convex if and only if it is
$r$-spherically supported and has \emph{nonempty interior}. This
equivalence was subsequently strengthened in
\cite[Theorem~4]{Pos2026}. Indeed, after it was proved in
\cite[Theorem~3]{Pos2026} that an $r$-spherically supported set is
either $r$-strongly convex or is contained in the boundary of a
nonempty closed $r$-strongly convex set, it was established in
\cite[Theorem~4]{Pos2026} that a nonempty closed set $S$ not reduced
to a singleton is $r$-strongly convex if and only if it is bounded,
has nonempty interior, and its boundary is $r$-spherically supported.

In \cite{Nacry2025}, Nacry and Thibault extended the notion of strong
convexity to variable radii by introducing the class of
$r(\cdot)$-strongly convex sets. For such sets, the supporting radius
at each boundary point $s$ is prescribed by the value $r(s)$ of a
positive function $r(\cdot)$ defined on the boundary. As in the
constant-radius setting, they established several geometric and
analytic characterizations of this class. This naturally raises the
following question: can the results of \cite{Pos2026} be extended to
the variable-radius setting by introducing a corresponding notion
of $r(\cdot)$-spherically supported sets?

The purpose of the present paper is to provide an affirmative answer
to this question, albeit with an essential modification involving an
\emph{induced} radius function. Given a positive function $r\colon\bdry S\longrightarrow ]0,+\infty[$, 
we say that $S$ is $r(\cdot)$-spherically supported if, for every
$s\in\bdry S$, there exists a unit proximal normal direction $\zeta$
to $S$ at $s$ such that \[S\subset\overline B\bigl(s-r(s)\zeta;r(s)\bigr).\] Unlike in the constant-radius case, the intersection of these
supporting balls need not be strongly convex with respect to
the prescribed function $r(\cdot)$. To address this
difficulty, we associate with $S$ a \emph{canonical} set\[ A_S:=
\bigcap_{(s,\zeta)\in G_S}
\overline B\bigl(s-r(s)\zeta;r(s)\bigr), \]
where $G_S$ records the boundary points of $S$ together with the unit
proximal normal vectors that are \emph{far realized} by the
corresponding $r(s)$-spheres, in the sense defined in
Section~\ref{Notation}. We also introduce an induced radius
function $\rho(\cdot)$ on the boundary of $A_S$. At each boundary point $x$,
the value $\rho(x)$ is the largest radius among the supporting balls
that are active at $x$, namely, those whose boundaries contain $x$.
This function reflects the local geometry of the canonical set more
accurately than the global radius \[
R:=\max\{r(s):s\;\hbox{is a boundary point of}\;S\},\] as illustrated in Example~\ref{example1}.

Our first main result establishes the fundamental properties of
$A_S$. Under the continuity of $r(\cdot)$, we prove that $S$, $A_S$,
and $G_S$ are compact and that $A_S$ is $R$-strongly convex. When $S$
is not a singleton, we obtain a representation of the normal cone
to $A_S$ in terms of the active supporting balls. We further prove
that $A_S$ is $\rho(\cdot)$-strongly convex and that $\rho(\cdot)$ is
upper semicontinuous. These properties lead to the principal geometric
characterization of the paper: under the continuity of $r(\cdot)$,
every nonempty closed $r(\cdot)$-spherically supported set $S$ satisfies
one of the following two conclusions:
\begin{itemize}
\item $S=A_S$ and $S$ is $r(\cdot)$-strongly convex;
\item $S$ is contained in the boundary of the nonempty, compact, and
      $\rho(\cdot)$-strongly convex set $A_S$.
\end{itemize}
The first assertion holds whenever $S$ has nonempty interior, whereas
the second holds whenever $S$ has empty interior. Thus, the dichotomy
known in the constant-radius setting persists for variable radii. In
the second alternative, however, the naturally induced function
$\rho(\cdot)$ replaces the prescribed function $r(\cdot)$. This
distinction is essential. We provide an explicit two-point example
showing that, in general, the induced radius function cannot be
replaced by the prescribed radius function at every point of $S$.
On the other hand, when $r(\cdot)$ is constant, we show that
$\rho(\cdot)$ coincides with $r$ and therefore recover
\cite[Theorem~3]{Pos2026}. We conclude with a characterization of
$r(\cdot)$-strongly convex sets that extends
\cite[Theorem~4]{Pos2026} to the variable-radius setting.

The paper is organized as follows. Section~\ref{Notation} introduces
the notation and the basic definitions from nonsmooth and convex
analysis used throughout the paper. In Section~\ref{MR}, we construct
the canonical set $A_S$ and the induced radius function $\rho(\cdot)$,
establish their main properties, and prove the geometric
characterization of $r(\cdot)$-spherically supported sets. The section
also contains examples illustrating the roles of the prescribed and
induced radii. Finally, Appendix~\ref{App} contains the proof of the
approximation result for extreme normal directions used in the main
characterization.

\section{Notation and basic definitions} \label{Notation} 

Throughout the paper, $\mathbb R^n$ is equipped with the Euclidean norm
$\|\cdot\|$ and the usual inner product $\langle\cdot,\cdot\rangle$. Its
open and closed unit balls are denoted by $B$ and $\overline B$,
respectively. We also denote by
$\mathbb S^{n-1}$ the unit sphere in $\mathbb R^n$. Given
$x\in\mathbb R^n$ and $r>0$, we write
\[
B(x;r):=x+rB
\;\;\text{and}\;\;
\overline B(x;r):=x+r\overline B.
\]
We denote by $\mathbb R_+:=[0,+\infty[$. For a set $S\subset\mathbb R^n$, its interior, closure and boundary are denoted by
$\inte S$, $\clo S$, and $\bdry S$, respectively, and its convex and conical hulls
by $\operatorname{co}S$ and $\operatorname{cone}S$, respectively. The
closed (resp.\ open) segment joining $x,y\in\mathbb R^n$ is denoted by
$[x,y]$ (resp.\ $]x,y[$). For a vector
$e\in\mathbb R^n$, its orthogonal complement is denoted by
\[
e^\perp:=\{x\in\mathbb R^n:\langle x,e\rangle=0\}.
\] 
For a nonempty set $S\subset\mathbb R^n$, the farthest distance from
$x\in\mathbb R^n$ to $S$ is defined by
\[
d_S^f(x):=\sup_{s\in S}\|x-s\|.
\]
The corresponding set of farthest points is denoted by $\far_S(x)$ and
given by
\[
\far_S(x):=
\bigl\{s\in S:\|x-s\|=d_S^f(x)\bigr\}.
\]

We next introduce the notions from nonsmooth and convex  analysis needed in the
sequel. For comprehensive treatments of this material, we refer to the
monographs \cite{clsw,mord,penot,rockwet,rockbook,thibault}. A convex body is a compact convex subset of $\mathbb R^n$ with
nonempty interior. Let $K\subset\mathbb R^n$ be a convex cone. The cone $K$ is said to be
\emph{pointed} if
\[
K\cap(-K)=\{0\},
\]
or, equivalently, if $K$ contains no line. A ray
$\mathbb R_+\zeta$, with nonzero $\zeta\in K$, is called an
\emph{extreme ray} of $K$ if
\[
\big[\zeta_1,\,\zeta_2\in K
\;\;\text{and}\;\;
\zeta_1+\zeta_2\in\mathbb R_+\zeta\big]\Longrightarrow
\zeta_1,\zeta_2\in\mathbb R_+\zeta.
\]
For a convex set $Q\subset\mathbb R^n$, a point $q\in Q$ is called an
\emph{extreme point} of $Q$ if
\[ \Big[q=(1-\lambda)q_1+\lambda q_2,\;\;q_1,q_2\in Q,\;\; \lambda\in]0,1[\Big] \Longrightarrow q_1=q_2=q.
\]
Let
$S\subset\mathbb R^n$ be a nonempty closed set and let $s\in S$. The
\emph{proximal normal cone} to $S$ at $s$, denoted by $N_S^P(s)$, is
defined by
\[
N_S^P(s):=
\left\{
\zeta\in\mathbb R^n: \exists \sigma\geq0\;\text{such that}\;\langle\zeta,x-s\rangle
\leq\sigma\|x-s\|^2
\;\;\text{for all }x\in S
\right\}.
\]
If $S$ is convex, the proximal normal cone $N_S^P(s)$ reduces to the
\emph{normal cone of convex analysis}, written $N_S(s)$ and given by
\[
N_S(s):=
\bigl\{
\zeta\in\mathbb R^n:
\langle\zeta,x-s\rangle\leq0
\;\;\text{for all }x\in S
\bigr\}.
\]
In this setting, $N_S(s)$ is nontrivial for every
$s\in\bdry S$. A point $s\in\bdry S$ is called \emph{regular} if
there exists a unit vector $\xi_s\in\mathbb R^n$ such that
\[
N_S(s)=\{\lambda\xi_s:\lambda\geq0\}.
\]
The corresponding supporting halfspace
\[
\bigl\{x\in\mathbb R^n:
\langle\xi_s,x-s\rangle\leq0\bigr\}
\]
is then called a \emph{regular supporting halfspace} of $S$. For $r>0$, a unit vector
$\zeta\in N_S^P(s)$ is said to be \emph{realized by an $r$-sphere} if
\[
B(s+r\zeta;r)\cap S=\emptyset,
\]
or, equivalently,
\[
\langle\zeta,x-s\rangle
\leq\frac{1}{2r}\|x-s\|^2
\;\;\text{for all }x\in S.
\]
The vector $\zeta$ is said to be \emph{far realized by an $r$-sphere}
if
\[
s\in\far_S(s-r\zeta),
\]
or, equivalently,
\[
\langle\zeta,x-s\rangle
\leq-\frac{1}{2r}\|x-s\|^2
\;\;\text{for all }x\in S,
\]
which is also equivalent to
\[
S\subset\overline B(s-r\zeta;r).
\]
For a nonempty closed set $S\subset\mathbb R^n$ and a function
$r\colon\bdry S\longrightarrow ]0,+\infty[$, we say that $S$ is
\emph{$r(\cdot)$-strongly convex} if, for every $s\in\bdry S$, every
unit vector $\zeta\in N_S^P(s)$ is far realized by an $r(s)$-sphere.
The characterizations established in \cite{Nacry2025} imply that every
$r(\cdot)$-strongly convex set is bounded and convex and, unless it is
reduced to a singleton, has nonempty interior. In comparison, $S$ is said to be
\emph{$r(\cdot)$-spherically supported} if, for every
$s\in\bdry S$, there exists a unit vector
$\zeta\in N_S^P(s)$ that is far realized by an $r(s)$-sphere. Thus,
$S$ is $r(\cdot)$-strongly convex if
\[
S\subset\overline B\bigl(s-r(s)\zeta;r(s)\bigr)
\]
for every boundary point $s$ and every unit proximal normal direction
$\zeta$ at $s$, whereas $S$ is $r(\cdot)$-spherically supported if the
same inclusion holds for at least one such direction at each boundary
point. When $r(\cdot)$ is constant, these notions reduce to
$r$-strong convexity and $r$-spherical support, respectively. We say that a nonempty closed set $S\subset\mathbb R^n$ satisfies the
\emph{exterior $\infty$-sphere condition} if, for every
$s\in\bdry S$, there exists a unit vector $\zeta\in N_S^P(s)$ that is
realized by an $r$-sphere for every $r>0$. Finally, for an extended real-valued function
$f\colon\mathbb R^n\to\mathbb R\cup\{+\infty\}$, its effective domain
and epigraph are defined, respectively, by
\[
\dom f:=\{x\in\mathbb R^n:f(x)<+\infty\}\;\;\text{and}\;\; \epi f:=\bigl\{(x,\alpha)\in\mathbb R^n\times\mathbb R:
f(x)\leq\alpha\bigr\}.
\]
The function $f$ is said to be \emph{proper} if $\dom f\neq\emptyset$,
and \emph{closed} if $\epi f$ is closed. If $f$ is proper and convex, its convex subdifferential at
$x\in\dom f$ is denoted by $\partial f(x)$ and defined by
\[
\partial f(x):=
\bigl\{
v\in\mathbb R^n:
f(y)\geq f(x)+\langle v,y-x\rangle
\;\;\text{for all }y\in\mathbb R^n
\bigr\}.
\]

\section{Main results}\label{MR} We begin by introducing the notation and the canonical
constructions that will be used throughout the sequel. For a nonempty closed set $S\subset\mathbb R^n$, a point
$s\in\bdry S$, and a function
$r\colon\bdry S\to]0,+\infty[$, we denote by
\[[N^P_S(s)]^f_{r(s)}:=\{\zeta\in N_S^P(s) : \zeta\;\hbox{is unit and far realized by an}\;r(s)\hbox{-sphere}\}.\]
We also set \[G_S:= \left\{ (s,\zeta): s\in \operatorname{bdry}S\;\hbox{and}\;\zeta\in [N_S^P(s)]^f_{r(s)} \right\}, \;\;A_S:= \bigcap_{(s,\zeta)\in G_S} \overline B(s-r(s)\zeta;r(s)),\] and  \[ R:=\sup_{s\in\bdry S}r(s).\]
For $\displaystyle x\in \operatorname{bdry}A_S\supset A_S\cap\bigcup_{(s,\zeta)\in G_S} \bdry \overline B(s-r(s)\zeta;r(s)) $, define \[G_S(x):= \left\{(s,\zeta)\in G_S : x\in \bdry \overline B(s-r(s)\zeta;r(s)) \right\}, \]
and set \begin{equation} \label{rho} \rho(x):= \sup_{(s,\zeta)\in G_S(x)} r(s)\leq R.\end{equation}
Finally, for every $(s,\zeta)\in G_S$, define
\[
\varphi_{s,\zeta}(x):=
\|x-s+r(s)\zeta\|^2-r(s)^2,\;\;\forall x\in\mathbb R^n,
\]
and set
\[
\Phi(x):=\sup_{(s,\zeta)\in G_S}\varphi_{s,\zeta}(x).
\]
Then
\[
A_S=\{x\in\mathbb R^n:\Phi(x)\le0\}.
\]

The following theorem is the cornerstone of our analysis. It establishes
the fundamental geometric and variational properties of the canonical set
$A_S$ and provides the foundation for the main characterization of
$r(\cdot)$-spherically supported sets.

\begin{theorem}[Structure of the canonical set $A_S$]\label{th0}
Let $S\subset\mathbb R^n$ be a nonempty closed set, and let $r\colon\bdry S\to]0,+\infty[$ be continuous. Assume that $S$ is $r(\cdot)$-spherically supported. Then $S$ is compact, $S\subset A_S$, $A_S$ is compact, $G_S$ is compact, and $A_S$ is $R$-strongly convex. Furthermore, if $S$ is not reduced to a singleton, then:
\begin{itemize}
\item $\inte A_S\not=\emptyset$.
\item $N_{A_S}(x) =\operatorname{cone}\left\{x-s+r(s)\zeta: (s,\zeta)\in G_S(x) \right\}$ for all $x\in\bdry A_S$.
\item $A_S$ is $\rho(\cdot)$-strongly convex, and the function $\rho(\cdot)$ defined in \eqref{rho} is upper semicontinuous. 
\end{itemize}
\end{theorem}

\begin{proof} As $S$ is nonempty and $r(\cdot)$-spherically supported, we deduce that  $G_S\neq\emptyset$. Hence, there exists $(s_0,\zeta_0)\in G_S$ such that \[S\subset A_S\subset \overline B(s_0-r(s_0)\zeta_0;r(s_0)). \]
Consequently, both $S$ and $A_S$ are bounded. Since $S$ is closed and $A_S$ is closed as an intersection of closed
balls, both sets are compact. Since $G_S\subset \bdry S\times  \mathbb S^{n-1}$, we have that $G_S$ is bounded. We claim that $G_S$ is closed. Let $(s_k,\zeta_k)\in G_S$ with $(s_k,\zeta_k)\f(s,\zeta)$. Since $r$ is continuous, passing to the limit in the far-realization inequalities shows that $\zeta\in [N_S^P(s)]^f_{r(s)}$. Hence, $(s,\zeta)\in G_S$. Thus, $G_S$ is closed, and consequently compact. We proceed to prove that $A_S$ is $R$-strongly convex. Since $r$ is continuous and finite-valued on the compact set $\bdry S$, the constant $R$ is well defined and finite. Clearly, the mapping
\[
(x,s,\zeta)\mapsto \|x-s+r(s)\zeta\|^2-r(s)^2
\]
is continuous on $\mathbb R^n\times G_S$. As $G_S$ is compact, it follows from Berge's maximum theorem that the function $
\Phi(\cdot)$ is continuous. Let $x\in\bdry A_S$. Since
\[
A_S=\{x\in\mathbb R^n:\Phi(x)\leq0\},
\]
and $\Phi$ is continuous, we have $\Phi(x)=0$. Moreover, the compactness
of $G_S$ ensures that the maximum defining $\Phi(x)$ is attained. Hence, there exists $(s,\zeta)\in G_S$ such that
\[
\varphi_{s,\zeta}(x)=\Phi(x)=0.
\]
Therefore,
\[
x\in \bdry \overline B(s-r(s)\zeta;r(s)).
\]
For $\nu:=\frac{x-s+r(s)\zeta}{r(s)}$, we have \[\nu\in N^P_{\overline B(s-r(s)\zeta;r(s))} (x)\;\;\hbox{and}\;\;A_S\subset \overline B(s-r(s)\zeta;r(s)).\]
Hence, $\nu\in N_{A_S}^P(x)$. Moreover, since $r(s)\le R$, the ball $\overline B(s-r(s)\zeta;r(s))$ is contained in the ball $
\overline B(x-R\nu;R)$. Indeed, the two balls are tangent at $x$ and have the same outward normal $\nu$,
with radii $r(s)\le R$. Therefore,
\[
A_S\subset \overline B(s-r(s)\zeta;r(s))\subset \overline B(x-R\nu;R).
\]
Thus, for every $x\in\operatorname{bdry}A_S$, there exists a unit vector
$\nu\in N_{A_S}^P(x)$ such that
\[
A_S\subset \overline B(x-R\nu;R).
\]
Hence, $A_S$ is $R$-spherically supported. Since $A_S$ is convex, being the intersection of closed balls, it follows from \cite[Corollary 1.3]{JCA2026} that $A_S$ is $R$-strongly convex. Assume now that $S$ is not reduced to a singleton. Since $S\subset A_S$ and every non-singleton $R$-strongly convex set has nonempty interior, we conclude that $\operatorname{int}A_S\neq\emptyset$. Then, there exists $y\in \inte A_S$. Since $y\in\inte A_S$, $A_S=\{x\in\R^n:\Phi(x)\leq0\}$, and
$G_S$ is compact, the supremum defining $\Phi(y)$ is attained. Hence,
$\Phi(y)<0$, and therefore the Slater condition holds. Therefore, by the normal cone formula for convex sublevel sets under
Slater's condition \cite[Theorem~23.7]{rockbook}, we have
\begin{equation} \label{normalA_S} N_{A_S}(x)=\mathbb R_+\partial\Phi(x)\;\;\text{for}\;x\in\bdry A_S.\end{equation}
Moreover, for every fixed $x\in \R^n$, the mapping \[(s,\zeta)\mapsto\varphi_{s,\zeta}(x)=\|x-s+r(s)\zeta\|^2-r(s)^2\]
is continuous on $G_S$. Together with the compactness of $G_S$ and the convexity and continuity of each function $\varphi_{s,\zeta}$, all the assumptions of \cite[Theorem~4.4.2]{HUL} are satisfied. Therefore,
\[
\partial\Phi(x)
=
\operatorname{co}
\left\{
\nabla\varphi_{s,\zeta}(x):
(s,\zeta)\in G_S(x)
\right\}.
\]
Since $
\nabla\varphi_{s,\zeta}(x)
=
2\left(x-s+r(s)\zeta\right)$, it follows from \eqref{normalA_S}  that
\[
N_{A_S}(x)
=
\operatorname{cone}
\left\{
x-s+r(s)\zeta:
(s,\zeta)\in G_S(x)
\right\}\;\text{for}\;x\in\operatorname{bdry}A_S.
\]
Now we prove that $A_S$ is $\rho(\cdot)$-strongly convex. Let $x\in\operatorname{bdry}A_S$ and let
$v\in N_{A_S}(x)$ with $\|v\|=1$. Then there exist $(s_i,\zeta_i)\in G_S(x)$ and $\lambda_i\ge0$, $i=1,\ldots,m$, such that
\[
v=
\sum_{i=1}^{m}
\lambda_i
\bigl(
x-s_i+r(s_i)\zeta_i
\bigr).
\]
For each $i$, since $A_S\subset\overline B(s_i-r(s_i)\zeta_i;r(s_i))$ and $x\in\bdry\overline B(s_i-r(s_i)\zeta_i;r(s_i))$, we have, for every $y\in A_S$,
\[
\left\langle
x-s_i+r(s_i)\zeta_i,
y-x
\right\rangle
\le
-\frac12
\|y-x\|^{2}.
\]
Multiplying by $\lambda_i$ and summing over $i=1,\ldots,m$, we obtain
\[
\langle v,y-x\rangle
\le
-\frac12
\left(
\sum_{i=1}^{m}\lambda_i
\right)
\|y-x\|^{2}.
\]
On the other hand,
\[
1=\|v\|
\le
\sum_{i=1}^{m}
\lambda_i
\|x-s_i+r(s_i)\zeta_i\|
=
\sum_{i=1}^{m}
\lambda_i r(s_i)
\le
\rho(x)
\sum_{i=1}^{m}\lambda_i.
\]
Therefore,
\[
\sum_{i=1}^{m}\lambda_i
\ge
\frac1{\rho(x)}.
\]
Consequently,
\[
\langle v,y-x\rangle
\le
-\frac1{2\rho(x)}
\|y-x\|^{2},\;\;
\forall y\in A_S.
\]
This proves that $A_S$ is
$\rho(\cdot)$-strongly convex. Finally, the upper semicontinuity of $\rho(\cdot)$ follows from the continuity of $r(\cdot)$, the compactness of $G_S$, and the closedness of the graph of the multifunction $x\mapsto G_S(x)$.  \end{proof}

\begin{remark} Under the hypotheses of Theorem~\ref{th0}, the preceding proof shows that the function $\Phi(\cdot)$ is continuous. Since $A_S=\{x\in\mathbb R^n:\Phi(x)\le0\}$, it follows that
\[\inte A_S=\{x\in\mathbb R^n:\Phi(x)<0\}=\bigcap_{(s,\zeta)\in G_S} B(s-r(s)\zeta;r(s)),\;\;\text{and}\]
\[\bdry A_S=\{x\in\mathbb R^n:\Phi(x)=0\}=A_S\cap \bigcup_{(s,\zeta)\in G_S} \bdry \overline B(s-r(s)\zeta;r(s)).\]
\end{remark}

The next example illustrates the geometric significance of the local radius function $\rho(\cdot)$ introduced in Theorem~\ref{th0}. In particular, it shows that the local radius may provide a strictly sharper description of the geometry than the global radius $R$, since one may have $\rho(x)<R$ at some boundary points.

\begin{figure}[h!]
\centering
\includegraphics[width=60mm]{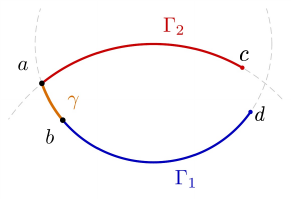}
\caption{\label{Fig1} The set $S$ of Example \ref{example1}}
\end{figure}

\begin{example}\label{example1}
Inspired by \cite[Figure~1]{Nacry2025}, consider the planar set
$S:=\Gamma_1\cup\Gamma_2$, shown in Fig.~\ref{Fig1}, where
$\Gamma_1$ and $\Gamma_2$ are two circular arcs meeting only at the
point $a$. Assume that the closed disk determined by each arc contains
$S$ and that its boundary circle has no other point of contact with the
opposite arc. Let $R_1$ and $R_2$ denote the radii of $\Gamma_1$ and
$\Gamma_2$, respectively, and suppose that
\[
R_2>R_1>0.
\]
Let $\gamma\subset\Gamma_1$ be a small subarc with $a$ as one
endpoint, and let $b$ denote the other endpoint of $\gamma$. Define the function $r\colon\bdry S\to]0,+\infty[$ by \[ r(s):= \begin{cases} R_1, & s\in\Gamma_1\cap\gamma^c,\\[1mm] R_2, & s\in\Gamma_2,\\[1mm] r_\gamma(s), & s\in\gamma\cap\{a\}^c, \end{cases} \]
where $r_\gamma\colon\gamma\to[R_1,R_2]$ is a continuous function such that
\[r_\gamma(b)=R_1\;\;\text{and}\;\;r_\gamma(a)=R_2.\] Then $r$ is continuous on $\bdry S$, and \[ R=\max_{s\in\bdry S}r(s)=R_2. \]
Clearly, $S$ is $r(\cdot)$-spherically supported. Since $S$ is not reduced to a singleton and $r(\cdot)$ is continuous,
Theorem~\ref{th0} yields that $A_S$ is nonempty, compact, and
$\rho(\cdot)$-strongly convex, and that $\rho(\cdot)$ is upper
semicontinuous. Let $d$ denote the other endpoint of $\Gamma_1$. By the construction
of the two arcs, for every
$x\in\Gamma_1\cap\gamma^c\cap\{d\}^c$, the only active supporting ball
is the disk determined by $\Gamma_1$. Therefore,
\[
\rho(x)=R_1<R.
\]
Hence, the local radius
function $\rho(\cdot)$ may provide a strictly sharper description of
the geometry than the global radius $R$.
\end{example}

We are now in a position to state the main result of this paper,
Theorem~\ref{th1}, which describes the geometric structure of
$r(\cdot)$-spherically supported sets through the canonical set $A_S$. Its proof hinges on the following extreme normal approximation lemma, whose proof is deferred to Appendix~\ref{App} for completeness.

\begin{lemma}[Extreme normal approximation]\label{lemnew}
Let $C\subset\mathbb R^n$ be a convex body, let $c\in\bdry C$, and let
$\zeta\in N_C(c)$ be a unit vector generating an extreme ray of
$N_C(c)$. Then there exist regular boundary points $c_k\in\bdry C$
converging to $c$ such that the unique unit outward normal vectors
$\zeta_k$ to $C$ at $c_k$ converge to $\zeta$.
\end{lemma}

\begin{theorem}[Geometric characterization of $r(\cdot)$-spherically supported sets]\label{th1} Let $S\subset\mathbb R^n$ be a nonempty closed set, and let $r\colon\bdry S\to]0,+\infty[$ be a continuous function. Assume that $S$ is $r(\cdot)$-spherically supported. Then one of the following alternatives holds:
\begin{enumerate}[$(i)$]
\item $S=A_S$ and is $r(\cdot)$-strongly convex.
\item $S\subset\bdry A_S$, where $A_S$ is nonempty, compact, and $\rho(\cdot)$-strongly convex.
\end{enumerate}
Here, $\rho(\cdot)$ is the upper semicontinuous function defined in \eqref{rho}. Moreover, assertion $(i)$ holds whenever $\inte S\neq\emptyset$, while assertion $(ii)$ holds whenever $\inte S=\emptyset$.
\end{theorem}

\begin{proof}  Let $S\subset\mathbb R^n$ be a nonempty closed set, and let $r\colon\bdry S\to]0,+\infty[$ be a continuous function such that $S$ is $r(\cdot)$-spherically supported. Then $S$ is compact, and for any $s\in\bdry S$, there exists a unit vector $\zeta\in N_S^P(s)$ such that \[\<\zeta,x-s\>\leq -\frac{1}{2r(s)}\|x-s\|^2,\;\;\forall x\in S.\]
Hence, for all $\delta>0$, \[\<\zeta,x-s\>\leq -\frac{1}{2r(s)}\|x-s\|^2\leq\frac{1}{2\delta}\|x-s\|^2,\;\;\forall x\in S.\]
This yields that $S$ satisfies the exterior $\infty$-sphere condition.\vspace{0.1cm}\\
{\underline{Case 1}}: $\inte S\not=\emptyset$.\vspace{0.1cm}\\
Then by \cite[Theorem 2]{Pos2026}, the set $S$ is convex. Hence, by \cite[Theorem 2.2.6]{schneider} combined with the definition of the $r(\cdot)$-spherically supported property, we have that \begin{equation*}A_S\supset S= \bigcap_{\substack{s\in\bdry S\\s\;\text{regular}\\\zeta_s\in N^P_S(s)\;\text{unit}}} \{x\in\R^n : \<\zeta_s,x-s\>\leq 0\} \supset \bigcap_{\substack{s\in\bdry S\\s\;\text{regular}\\ \zeta_s\in [N^P_S(s)]^f_{r(s)}}} \widebar{B}(s-r(s)\zeta_s;r(s)) \supset A_S \end{equation*}
This yields that $S=A_S$. We proceed to prove that $S$ is $r(\cdot)$-strongly convex. Fix
$s\in\bdry S$ and let $\zeta\in N_S(s)$ be a unit vector generating
an extreme ray of $N_S(s)$. By Lemma~\ref{lemnew}, there exist regular
boundary points $s_k\in\bdry S$ converging to $s$ such that the unique
unit outward normal vectors $\zeta_k$ to $S$ at $s_k$ converge to
$\zeta$. Since $S$ is $r(\cdot)$-spherically supported and $s_k$ is
regular, $\zeta_k$ is far realized by an $r(s_k)$-sphere. Therefore,
\[
\langle\zeta_k,x-s_k\rangle
\leq-\frac{1}{2r(s_k)}\|x-s_k\|^2
\;\;\text{for all }x\in S.
\]
Since $r(\cdot)$ is continuous, $r(s_k)\to r(s)$. Passing to the limit
in the preceding inequality yields
\begin{equation}\label{eq:extreme-normal}
\langle\zeta,x-s\rangle
\leq-\frac{1}{2r(s)}\|x-s\|^2
\;\;\text{for all }x\in S.
\end{equation}
Thus, every unit vector generating an extreme ray of $N_S(s)$ is far
realized by an $r(s)$-sphere. Now let $\zeta\in N_S(s)$ be an arbitrary unit vector. Since $S$ is convex and has nonempty interior, the closed convex cone
$N_S(s)$ is pointed. Hence, by
\cite[Corollary~18.5.2]{rockbook}, it is generated by its extreme
rays. Consequently, there exist unit vectors
$\zeta_1,\ldots,\zeta_m$, each generating an extreme ray of $N_S(s)$,
and numbers $\lambda_1,\ldots,\lambda_m\geq0$ such that
\[
\zeta=\sum_{i=1}^m\lambda_i\zeta_i.
\]
Since $\zeta$ and the vectors $\zeta_i$ are unit, the triangle
inequality gives
\[
1=\|\zeta\|
\leq\sum_{i=1}^m\lambda_i.
\]
Applying \eqref{eq:extreme-normal} to each $\zeta_i$, we obtain, for
every $x\in S$,
\[
\langle\zeta,x-s\rangle=\sum_{i=1}^m\lambda_i\langle\zeta_i,x-s\rangle \leq-\frac{1}{2r(s)}
\left(\sum_{i=1}^m\lambda_i\right)\|x-s\|^2 \leq-\frac{1}{2r(s)}\|x-s\|^2.
\]
Equivalently,
\[
S\subset\overline B\bigl(s-r(s)\zeta;r(s)\bigr).
\]
Since $s\in\bdry S$ and the unit vector $\zeta\in N_S(s)$ were
arbitrary, $S$ is $r(\cdot)$-strongly convex.\vspace{0.1cm}\\
{\underline{Case 2}}: $\inte S=\emptyset$.\vspace{0.1cm}\\
Let $s\in\bdry S=S\subset A_S$. As $S$ is $r(\cdot)$-spherically supported, there exists $\zeta_s\in[N^P_S(s)]^f_{r(s)}$. We have  $(s,\zeta_s)\in G_S$ and \[s\in A_S\cap \bdry \overline B(s-r(s)\zeta_s;r(s))\subset \bdry A_S.\]
This yields that $S\subset \bdry A_S$. If $S=\{s_0\}$ is a singleton, then clearly $A_S=S$ and $A_S$ is $\rho(\cdot)$-strongly convex, with $\rho(\cdot)\equiv r(s_0)$. Otherwise, Theorem~\ref{th0} yields that $A_S$ is nonempty, compact, and $\rho(\cdot)$-strongly convex, with $\rho(\cdot)$ upper semicontinuous. \end{proof}

The following example shows that, in general, the induced radius
function $\rho(\cdot)$ appearing in Theorem~\ref{th1}$(ii)$ cannot be
replaced by the prescribed radius function $r(\cdot)$ at every point of
$S$. Consequently, unlike the constant-radius case, this theorem does
not directly yield an equivalent characterization in terms of the
original radius function $r(\cdot)$.

\begin{figure}[h!]
\centering
\includegraphics[width=80mm]{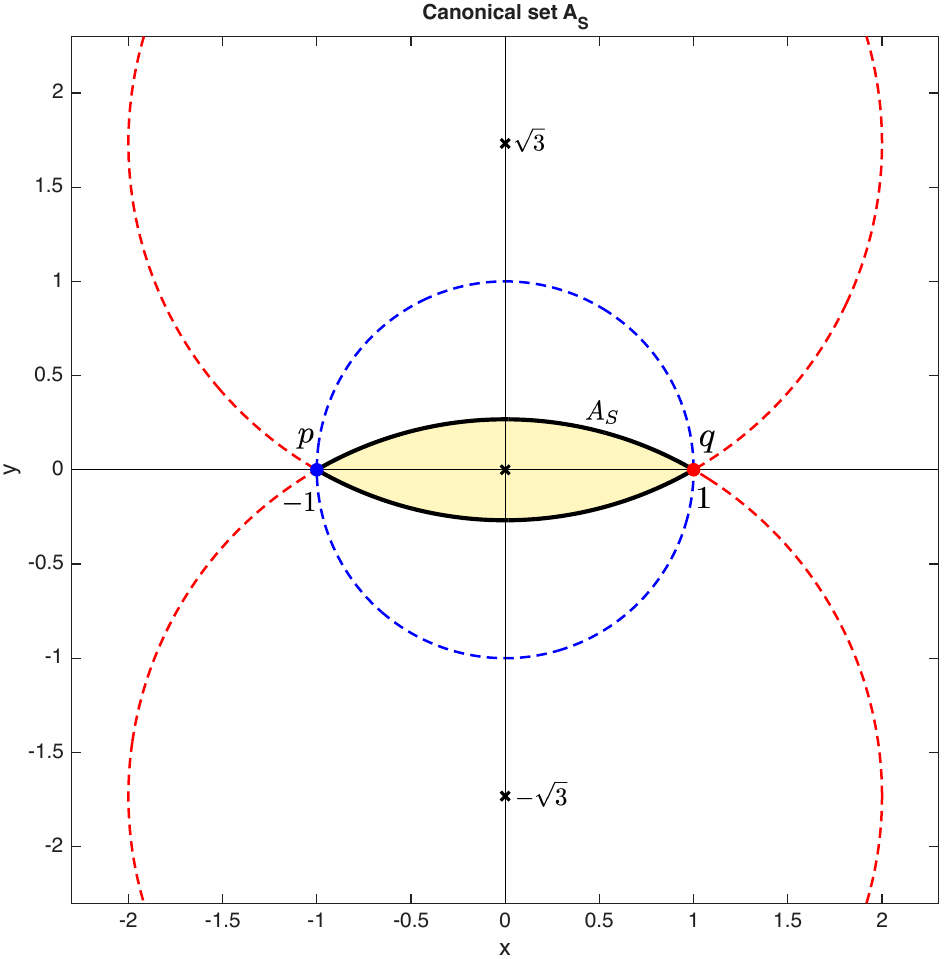}
\caption{\label{Fig2} The canonical set $A_S$ of Example \ref{example2}}
\end{figure}

\begin{example}\label{example2}
Consider the set $S:=\{p,q\}\subset\mathbb R^2$, where $p=(-1,0)$ and $q=(1,0)$, and define $r(p)=1$ and $r(q)=2$. Clearly, $r(\cdot)$ is continuous, and $S$ is $r(\cdot)$-spherically supported. Moreover, the active sets at $p$ and $q$ are given by
\[
G_S(p)=\left\{\left(p,(-1,0)\right),\left(q,\left(\frac12,\frac{\sqrt3}{2}\right)\right),\left(q,\left(\frac12,-\frac{\sqrt3}{2}\right)\right)\right\},\]
while \[G_S(q)=\left\{\left(p,(-1,0)\right)\right\}\cup\left\{(q,\zeta):\zeta\in\mathbb S^1,\ \langle \zeta,(1,0)\rangle\ge \frac12
\right\}.\]
A direct computation shows that
\[
A_S
=
\overline B((0,\sqrt3);2)
\cap
\overline B((0,-\sqrt3);2).
\]
Hence, $A_S$ is a lens-shaped set bounded by two circular arcs of radius
$2$, as shown in Fig.~\ref{Fig2}. Consequently,
\[
\rho(p)=2>1=r(p)
\;\;\text{and}\;\;
\rho(q)=2=r(q).
\]
Moreover, the prescribed radius $r(p)=1$ cannot be used at $p$ to
establish strong convexity of $A_S$ in every normal direction. Indeed, the vector
\[
\nu:=\left(-\frac12,-\frac{\sqrt3}{2}\right)
\]
is a unit outward normal to $A_S$ at $p$, whereas
\[
\|q-(p-\nu)\|=\sqrt3>1.
\]
Therefore, the ball $\overline B(p-\nu;1)$ does not contain
$q\in A_S$. Thus, the induced radius function $\rho(\cdot)$ cannot, in
general, be replaced by the prescribed radius function $r(\cdot)$ in
the second alternative of Theorem~\ref{th1}.
\end{example}

The following corollary recovers the constant-radius characterization established in \cite[Theorem~3]{Pos2026}.

\begin{corollary}[{\cite[Theorem 3]{Pos2026}}]Let $S\subset\R^n$ be a nonempty closed set, and let $r>0$. The following assertions are equivalent:
\begin{enumerate}[$(i)$]
\item $S$ is $r$-spherically supported.
\item $S$ is $r$-strongly convex, or there exists a nonempty, closed, and $r$-strongly convex set $A\subset\R^n$ such that $S\subset\bdry A$.
\end{enumerate}
Moreover, if $(i)$ holds, then $S$ is $r$-strongly convex whenever $\inte S\neq\emptyset$, and $S$ is contained in the boundary of a nonempty, closed, and $r$-strongly convex set whenever $\inte S=\emptyset$.
\end{corollary}
\begin{proof} The implication $(ii)\im(i)$ is immediate. For $(i)\Longrightarrow(ii)$ and the additional conclusions, it suffices to apply Theorem~\ref{th1} and observe that
$\rho(\cdot)\equiv r$. \end{proof}

The next corollary provides a characterization of $r(\cdot)$-strongly convex sets and extends \cite[Theorem~4]{Pos2026} to the variable-radius setting.

\begin{corollary}[Characterization of $r(\cdot)$-strongly convex sets] Let $S\subset\mathbb R^n$ be a nonempty closed set not reduced to a
singleton, and let $r\colon\bdry S\to]0,+\infty[$ be a continuous function. Then $S$ is $r(\cdot)$-strongly convex if and only if
$S$ is bounded, $\inte S\neq\emptyset$, and $\bdry S$ is
$r(\cdot)$-spherically supported.
\end{corollary}

\begin{proof} Assume first that $S$ is $r(\cdot)$-strongly convex. Then $S$ is bounded and, since it is not reduced to a singleton, has
nonempty interior. Moreover, the inclusion
\[
[N_S^P(\cdot)]_{r(\cdot)}^f
\subset
[N_{\bdry S}^P(\cdot)]_{r(\cdot)}^f
\]
shows that $\bdry S$ is $r(\cdot)$-spherically supported. Conversely, assume that $S$ is bounded, $\bdry S$ is $r(\cdot)$-spherically supported, and $\inte S\neq\emptyset$. We first show that $S$ is $r(\cdot)$-spherically supported. Fix
$s\in\bdry S$. Since $\bdry S$ is $r(\cdot)$-spherically supported,
there exists a unit vector
$\zeta\in N_{\bdry S}^P(s)$ such that
\begin{equation}\label{eq:boundary-supporting-ball}
\bdry S\subset
\overline B\bigl(s-r(s)\zeta;r(s)\bigr).
\end{equation}
We claim that
\begin{equation}\label{eq:S-convex-hull-boundary}
S\subset\operatorname{co}(\bdry S).
\end{equation}
Indeed, let $x\in S$. If $x\in\bdry S$, the claim is immediate. Suppose
that $x\in\inte S$, and fix a unit vector $u\in\mathbb R^n$. Since $S$
is closed and bounded, there exist $t_+,t_->0$ such that
\[
x_+:=x+t_+u\in\bdry S,\;\;
x_-:=x-t_-u\in\bdry S,
\]
and
\[
x=\frac{t_-}{t_++t_-}x_+
  +\frac{t_+}{t_++t_-}x_-.
\]
Thus, $x\in[x_-,x_+]\subset\operatorname{co}(\bdry S)$, which proves
\eqref{eq:S-convex-hull-boundary}. The ball appearing in \eqref{eq:boundary-supporting-ball} is convex.
Therefore,
\[
\operatorname{co}(\bdry S)
\subset
\overline B\bigl(s-r(s)\zeta;r(s)\bigr).
\]
Combining this inclusion with \eqref{eq:S-convex-hull-boundary}, we
obtain
\[
S\subset
\overline B\bigl(s-r(s)\zeta;r(s)\bigr).
\]
Equivalently,
\[
\langle\zeta,x-s\rangle
\leq-\frac{1}{2r(s)}\|x-s\|^2\;\;\text{for all }x\in S.
\]
In particular, $\zeta\in N_S^P(s)$ and is far realized by an
$r(s)$-sphere. Since $s\in\bdry S$ was arbitrary, $S$ is
$r(\cdot)$-spherically supported. As $\inte S\neq\emptyset$ and
$r(\cdot)$ is continuous, Theorem~\ref{th1}$(i)$ implies that $S$ is
$r(\cdot)$-strongly convex.
\end{proof}

\section*{Declarations}

\subsection*{Funding}
No funding was received for conducting this study.

\subsection*{Conflict of interest}
The author has no relevant financial or non-financial interests to
disclose.

\subsection*{Ethics approval and consent to participate}
Not applicable.

\subsection*{Clinical trial registration}
Not applicable.

\subsection*{Consent for publication}
Not applicable.

\subsection*{Data availability}
No datasets were generated or analyzed during the current study.

\subsection*{Materials availability}
Not applicable.

\subsection*{Code availability}
Not applicable.

\subsection*{Author contribution}
Chadi Nour conceived the study, developed the mathematical results and
proofs, prepared the figures, and wrote and revised the manuscript.

\begin{appendices}

\section{Proof of Lemma \ref{lemnew}}\label{App}

Fix $a\in\inte C$ and set
\[
e:=\frac{a-c}{\|a-c\|}.
\]
We first observe that
\begin{equation}\label{eq:strict-normal-direction}
\langle v,e\rangle<0
\;\;\text{for every nonzero } v\in N_C(c).
\end{equation}
Indeed, since $a\in\inte C$, for every nonzero $v\in N_C(c)$ there
exists $\varepsilon>0$ such that
\[
a+\varepsilon\frac{v}{\|v\|}\in C.
\]
It follows from the definition of $N_C(c)$ that
\[
\left\langle
v,a+\varepsilon\frac{v}{\|v\|}-c
\right\rangle\leq0,
\]
and hence
\[
\langle v,a-c\rangle\leq-\varepsilon\|v\|<0,
\]
which proves \eqref{eq:strict-normal-direction}. Write
\[
\mathbb R^n=e^\perp\oplus\mathbb Re.
\]
Let $P$ be the orthogonal projection of $C-c$ onto $e^\perp$. Since
$a-c$ is a positive multiple of $e$ and $a\in\inte C$, we have
$0\in\inte P$, where the interior is taken in $e^\perp$. Define the
extended real-valued function $f\colon e^\perp\to
\mathbb R\cup\{+\infty\}$ by
\[
f(y):=
\min\{t\in\mathbb R:c+y+te\in C\}
\quad\text{if }y\in P,
\]
and $f(y):=+\infty$ if $y\notin P$. The minimum exists because $C$ is
compact. Identifying $e^\perp\times\mathbb R$ with $\mathbb R^n$ through the
mapping $(y,t)\mapsto y+te$, we have
\[
\epi f=(C-c)+\mathbb R_+e,
\]
and therefore $f$ is a closed proper convex function. In view of
\eqref{eq:strict-normal-direction}, the point $c$ is the first point at which the oriented line
$c+\mathbb Re$ enters $C$ in the direction $e$. Thus, $f(0)=0$. For every $y\in\inte P$ and every $p\in e^\perp$, the definition of the
convex subdifferential gives
\[
p\in\partial f(y)\Longleftrightarrow
f(z)\geq f(y)+\langle p,z-y\rangle
\;\;\text{for every }z\in e^\perp.
\]
Consequently,
\begin{equation}\label{eq:normal-subdifferential}
N_C\bigl(c+y+f(y)e\bigr)
=
\bigl\{\lambda(p-e):\lambda\geq0,\;
p\in\partial f(y)\bigr\},\;\;\forall y\in\inte P.
\end{equation}
Indeed, if $p\in\partial f(y)$ and $c+z+te\in C$, then $t\geq f(z)$
and
\[\langle p-e,(c+z+te)-(c+y+f(y)e)\rangle =\langle p,z-y\rangle-t+f(y)
\leq f(z)-t
\leq0.
\]
Thus, $p-e\in N_C(c+y+f(y)e)$. Conversely, let $v\in N_C\bigl(c+y+f(y)e\bigr)$ be a nonzero vector. Since $C-c$ is convex with nonempty interior, the orthogonal
projection of $\inte(C-c)$ onto $e^\perp$ is $\inte P$. Hence, there
exists $t>f(y)$ such that $c+y+te\in\inte C$. The normal inequality then gives
\[
(t-f(y))\langle v,e\rangle\leq0,
\]
and therefore $\langle v,e\rangle\leq0$. If
$\langle v,e\rangle=0$, write $v=q+\alpha e$, where
$q\in e^\perp$. Then $\alpha=0$, and applying the normal inequality to a point of $C$ whose projection
is an arbitrary $z\in P$ gives
\[
\langle q,z-y\rangle\leq0\;\;\text{for every }z\in P.
\]
Since $y\in\inte P$, this implies $q=0$, contradicting $v\neq0$.
Consequently, $\langle v,e\rangle<0$. After multiplying $v$ by
a positive scalar, we may write $v=p-e$ for some $p\in e^\perp$.
Applying the normal inequality to the points $c+z+f(z)e$, for
$z\in P$, gives
\[
f(z)\geq f(y)+\langle p,z-y\rangle.
\]
For $z\notin P$, the same inequality holds automatically since
$f(z)=+\infty$. Thus, $p\in\partial f(y)$, and hence
\eqref{eq:normal-subdifferential} follows. In particular,
\begin{equation}\label{eq:normal-subdifferential-at-zero}
N_C(c)
=
\bigl\{\lambda(p-e):\lambda\geq0,\;
p\in\partial f(0)\bigr\}.
\end{equation}
The affine mapping $p\mapsto p-e$ identifies $\partial f(0)$ with
the section
\[
N_C(c)\cap
\bigl\{v\in\mathbb R^n:\langle v,e\rangle=-1\bigr\}.
\]
It follows that the ray generated by $p-e$ is an extreme ray of
$N_C(c)$ if and only if $p$ is an extreme point of $\partial f(0)$.
Therefore, since $\zeta$ generates an extreme ray of $N_C(c)$, there
exists an extreme point $p$ of $\partial f(0)$ such that
\begin{equation}\label{eq:zeta-p}
\zeta=\frac{p-e}{\|p-e\|}.
\end{equation}
Let $D$ denote the set of points in $\inte P$ at which $f$ is
differentiable, and let
\[
L(0):=
\left\{
\lim_{k\to\infty}\nabla f(y_k):
y_k\in D,\ y_k\to0
\right\}.
\]
By \cite[Theorem~25.5]{rockbook}, the set $D$ is dense in $\inte P$.
Moreover, since $0\in\inte(\dom f)$,
\cite[Theorem~25.6]{rockbook} yields
\begin{equation}\label{eq:subdifferential-gradient-limits}
\partial f(0)=\clo\bigl(\operatorname{co} L(0)\bigr).
\end{equation}
Since $f$ is finite and convex on a neighborhood of $0$, it is locally
Lipschitz there. Hence, its gradients at differentiability points in a
sufficiently small neighborhood of $0$ are uniformly bounded. It
follows that $L(0)$ is bounded. The set $L(0)$ is also closed. Indeed,
let $p_j\in L(0)$ and suppose that $p_j\to p_0$. By the definition of
$L(0)$, for every $j$ there exists $y_j\in D$ such that
\[
\|y_j\|<\frac1j
\;\;\text{and}\;\;
\|\nabla f(y_j)-p_j\|<\frac1j.
\]
Then $y_j\to0$ and $\nabla f(y_j)\to p_0$, so $p_0\in L(0)$. Thus,
$L(0)$ is compact, and \eqref{eq:subdifferential-gradient-limits}
becomes $\partial f(0)=\operatorname{co} L(0)$. Since $p$ is an extreme point of $\partial f(0)$, we necessarily have
$p\in L(0)$. Consequently, there exists a sequence $y_k\in D$ such
that
\begin{equation}\label{eq:gradient-convergence}
y_k\f 0
\;\;\text{and}\;\;
\nabla f(y_k)\f p.
\end{equation}
Set $c_k:=c+y_k+f(y_k)e$. Then $c_k\in\bdry C$. Since $f$ is continuous in a neighborhood of
$0$, it follows from $y_k\to0$ and $f(0)=0$ that $c_k\f c$. As $f$ is differentiable at $y_k$, we have
\[
\partial f(y_k)=\{\nabla f(y_k)\}.
\]
Thus, by \eqref{eq:normal-subdifferential},
\[
N_C(c_k)
=
\bigl\{\lambda(\nabla f(y_k)-e):\lambda\geq0\bigr\}.
\]
Hence, $c_k$ is a regular boundary point of $C$, and its unique unit
outward normal vector is
\[
\zeta_k:=
\frac{\nabla f(y_k)-e}
{\|\nabla f(y_k)-e\|}.
\]
Finally, \eqref{eq:zeta-p} and
\eqref{eq:gradient-convergence} imply that
\[
\zeta_k\f
\frac{p-e}{\|p-e\|}
=\zeta.
\]
The proof is complete. \hfill\qedsymbol
\end{appendices}



\begin{thebibliography}{9}

\bibitem{Azagra} Azagra, D., Drake, M., Hajlasz, P.: $C^2$-Lusin Approximation of Strongly Convex Bodies, J. Convex Anal. 32(1), 91--106 (2025)

\bibitem{clsw} Clarke, F. H., Ledyaev, Yu., Stern, R. J., Wolenski, P. R.: Nonsmooth Analysis and Control Theory, Graduate Texts in Mathematics 178, Springer, New York (1998)

\bibitem{frankowska} Frankowska, H., Olech, C.: $R$-convexity of the integral of set-valued functions, in: Contributions to Analysis and Geometry, Johns Hopkins Univ. Press, Baltimore, Md., 117--129 (1981)

\bibitem{HUL} Hiriart-Urruty, J.-B., Lemaréchal, C.: Fundamentals of Convex Analysis,  Springer-Verlag (2001)

\bibitem{mord} Mordukhovich, B. S.: Variational Analysis and Generalized Differentiation I: Basic Theory, Springer, Berlin (2006)

\bibitem{Nacry2025} Nacry, F., Thibault, L.: Strongly convex sets with variable radii, Math. Control Relat. Fields, 14(4), 1530--1559 (2024)

\bibitem{Pos2026}  Nour, C.: Spherically supported property and the exterior sphere condition with infinite radius, Positivity, to appear

\bibitem{JCA2026} Nour, C., Takche, J.: New characterization of strong convexity, J. Convex Anal. 34(1), 129--141 (2027)

\bibitem{penot} Penot, J.-P.: Calculus Without Derivatives, Graduate Texts in Mathematics 266, Springer, New York (2013)

\bibitem{rockbook} Rockafellar, R. T.: Convex Analysis, Princeton University Press (1970)

\bibitem{rockwet} Rockafellar, R. T., Wets, R. J.-B.: Variational Analysis, Grundlehren der Mathematischen Wissenschaften 317, Springer, Berlin (1998)

\bibitem{schneider} Schneider, R.: Convex Bodies: The Brunn–Minkowski Theory, 2nd ed., Cambridge Univ. Press (2013)


\bibitem{thibault} Thibault, L.: Unilateral Variational Analysis in Banach Spaces, World Scientific (2023)


\bibitem{weber} Weber, A., Reissig, G.: Local characterization of strongly convex sets, J. Math. Anal. Appl. 400, 743--750 (2013)

\end{thebibliography}

\end{document}